\documentclass [12pt]{article}
\usepackage{latexsym,amsthm,amsfonts,amsmath,url}

\DeclareMathOperator{\diam}{diam}

\newtheorem{theo}{Theorem}[section]

\newtheorem{coro}[theo]{Corollary}
\newtheorem{conj}[theo]{Conjecture}

\newcommand\Z{\ensuremath{\mathbb{Z}}}
\renewcommand\P{\ensuremath{\mathbb{P}}}

\title{A Resistance Bound via an Isoperimetric Inequality}
\author{ Itai Benjamini \and Gady Kozma}

\date{}
\begin{document}
\maketitle

\begin{abstract}
An isoperimetric  upper bound on the resistance is given. As a
corollary we resolve two problems, regarding mean commute time on
finite graphs and resistance on percolation clusters. Further
conjectures are presented.
\end{abstract}

\section{Introduction}
It is natural and useful to interpret a graph as an electrical
network by identifying the edges of the graph with one Ohm
resistors. Then the effective resistance between vertices and sets
admits a probabilistic and potential theoretic meaning, and is of
interest, see for instance \cite{DS}, \cite{P}, \cite{LP} for the
standard background and definitions. In the next section we
present a useful upper bound for the effective resistance between
two vertices in the graph in terms of an isoperimetric quantity
for connected sets containing one of these two vertices. We
suspect that this bound can't be truly new, still in
section~\ref{app} we bring two new applications.

\section{A resistance bound}\label{2}
For a subset $A$ of a graph $G$ we denote by $\partial A$ the \emph{external}
boundary, i.e.\ the vertices of $G\setminus A$ with neighbors in $A$. As usual, $\left\lfloor s\right\rfloor $ denotes the
largest integer $\leq s$ and $\left\lceil s\right\rceil $ denotes the
smallest integer $\geq s$. $\log$ denotes the logarithm to base $2$.

\begin{theo}
\label{main}
 Let $G$ be a finite graph. Let $w$ and
$u$ be vertices of $G$. Let $R_{w,u}$ be the electric resistance
between $w$ and $u$. Then\begin{equation}
R_{w,u}\leq C(L_{w}+L_{u}) , \quad L_{v}:=\sum _{n=1}^{\left\lfloor \log |G|\right\rfloor }\max _{\substack{ v\in A\\
 A\textrm{ connected}\\
 |G|2^{-(n+1)}<|A|\leq |G|2^{-n}}
}\frac{|A|}{|\partial A|^{2}}+\frac{1}{|\partial A|}
\label{eq:defrn}\end{equation}

\end{theo}

The $1/|\partial A|$ summands are relevant, of course, only in graphs with very
high connectivity. For example, let $G$ be a graph with $2n^2+2n+2$ vertices,
arranged as follows: $G=A_1\cup \dotsb \cup A_5$, $|A_1|=|A_5|=1$,
$|A_2|=|A_4|=n$ and $|A_3|=n^2$, and every vertex of $A_i$ is connected to
every vertex of $A_{i+1}$ and to every vertex of $A_{i-1}$. It is easy to see
that for this graph the resistance between $A_1$ and $A_5$ is $O(n^{-1})$, but
\[\sum \max \frac{|A|}{|\partial A|^2}=O(n^{-2})\quad .\]

To understand the conditions better it is worthwhile to examine
the following silly example: Let $G$ contain two components of
equal size with $w$ in one and $u$ in the other. In this case the
resistance is infinite. To get $\infty $ on the right hand side of
(\ref{eq:defrn}) you need to have a set $A$ with $|\partial A|=0$
and the only such set is the complete half. Therefore replacing
$|A|\leq |G|2^{-n}$ with $|A|\leq |G|2^{-n}-1$ would render the
lemma incorrect with any constant $C$.

\begin{proof}
We may assume that $w$ and $u$ are in the same component of $G$
since otherwise both sides of (\ref{eq:defrn}) are infinite. We
may assume that $G$ has no other components and pay a price of $2$
in the constant in (\ref{eq:defrn}). Connect a battery to $w$ and
$u$ so that the voltage $V$ satisfies $V(w)=1$ and $V(u)=0$.
Denote the electric current by $I$. Denote by $A_{m}$ a set of the
$m$ vertices with lowest voltage. If there are a number of
possibilities (because some vertices have equal voltage), we
choose $A_{m}$ to be connected, which can always be done due to
the maximum principle. Let $\theta (m)=\max _{v\in A_{m}}V(v)$.
Denote
\[ r_{n}:=\min \{|\partial A|:u\in A\textrm{, }A\textrm{ connected
and }|G|2^{-(n+1)}<|A|\leq |G|2^{-n}\}\quad .\] Let
$|G|2^{-(n+1)}<m\leq |G|2^{-n}$. As already remarked, the set
$A_{m}$ is connected due to the maximum principal, contains $u$
and therefore $|\partial A_{m}|\geq r_{n}$. For every vertex of $\partial A_m$
examine the sum of the currents through all edges connecting it to $A_m$ which
we will call for simplicity the current through the vertex. Since the average
current going through every vertex of $\partial A_{m}$ is $I/|\partial A_m|$
we get that for at least $\frac{1}{2}r_{n}$ vertices the current
through each is \[ \leq 2\frac{I}{|\partial A_m|}\leq
\frac{2I}{r_{n}}\quad .\]
Examine one such vertex $v$, and take an edge connecting $v$ to $A_m$. The
current through this edge, which is the difference of voltages, is $\leq
2I/r_n$. We get at least
$\frac{1}{2}r_{n}$ vertices $v$ where $V(v)\leq \theta
(m)+\frac{2I}{r_{n}}$. This gives that $\theta \big(m+\left\lceil\frac{1}{2}r_{n}\right\rceil \big)\leq \theta
(m)+\frac{2I}{r_{n}}$. We apply this to a series of $m$'s from
$|G|2^{-(n+1)}$ to $|G|2^{-n}$ and get
\[ \theta \left(\left\lfloor
|G|2^{-n}\right\rfloor +1\right)-\theta \left(\left\lfloor
|G|2^{-(n+1)}\right\rfloor +1\right)\leq
\frac{2I}{r_{n}}\left(\frac{2\left\lceil
|G|2^{-(n+1)}\right\rceil }{r_{n}}+1\right)\]
which we sum over $n$ and get \[
\theta (\left\lfloor |G|/2\right\rfloor +1)-\theta (1)\leq
CIL_{u}\quad .\] Since $\theta (1)=0$ we are done with the
neighborhood of $u$. An identical calculation around $w$ will
show\[ 1-\theta '(\left\lfloor |G|/2\right\rfloor +1)\leq
CIL_{w}\] where $\theta '(m)=\min _{v\in A_{m}'}V(m)$ and $A_{m}'$
is a set of $m$ vertices with maximal voltage which is a connected
neighborhood of $w$. Since $A_{\left\lfloor |G|/2\right\rfloor
+1}$ and $A_{\left\lfloor |G|/2\right\rfloor +1}'$ intersect we
get $\theta (\left\lfloor |G|/2\right\rfloor +1)\geq \theta
'(\left\lfloor |G|/2\right\rfloor +1)$ and then \[ 1\leq
CI(L_{u}+L_{w})\] which finishes the proof.\end{proof}
\medskip
\section{Applications}\label{app}
\subsection{Mean commute time}

Let $\tau^* = \max_{v,u}(E_vT_u + E_uT_v)$. Where $E_vT_u$ is the
expected hitting time for a random walk starting at $v$ to hit
$u$. By hitting time we mean here continuous hitting time, i.e.\ one puts on
every edge an alarm clock with the ringing time distributed like an exponential
variable and then move from a vertex along the first edge that rings, at the
time it rings. $\tau^*$ is the maximal mean commute time.

Open problem 20 in chapter 6 of \cite{AF} asserts the following:

\begin{quote}
Show that for real $1/2 < \gamma < 1$ and $\delta > 0$,
there exists a constant $K_{\gamma, \delta}$ with the following
property. Let $G$ be a regular $n$-vertex graph such that, for any
subset $A$ of vertices with $|A| \leq n/2$, there exist at least
$\delta|A|^{\gamma}$ edges between $A$ and $A^c$. Then $\tau^*
\leq K_{\gamma, \delta} n$.
\end{quote}

Theorem \ref{main} allows to prove this under the assumption that the graph
degree is bounded, in which case there is no difference between an
isoperimetric condition phrased in terms of the number of edges (as in the
problem) or in terms of the number of vertices (as in theorem~\ref{main}). In
this case we use the fact that the mean commute time between any two vertices
$u$, $w$ for simple random walk on a connected graphs equals $R_{u,w} |G|$, see
chapter 4 of \cite{AF}. The isoperimetric condition gives in theorem~\ref{main}
a bounded sum and the answer is positive.

When one removes the assumption of bounded degree, the answer is negative. For
example, take a graph with $n$ vertices arranged in a circle such that each two
neighbors are connected by $\lceil n^{2/3}\rceil$ edges. Then clearly the
assumptions hold (with $\gamma = \frac{2}{3}$ and $\delta=1$) but the
conclusion fails as $\tau^* >cn^{4/3}$. It is not difficult to construct such
an example with no multiple edges.

% Take the group Z_n X Z_{n^0.4} with the generator set {1,-1} X Z_n^0.4.
% Then the Cayley graph (which has degree n^0.4) satisfies the isoperimetric
% inequality but the resistance is still n^0.2. To see that it satisfies the
% isoperimetric inequality, define X to be the set of all n's such that
% the cut has size > 0.5n^0.4. If X is nonempty then the boundary is
% >0.25n^0.8. If X is empty then the boundary is > 0.5n^0.4|Omega|. If X is
% full then |Omega|>0.5n^1.4 i.e. 0.5|G|.
%
% In this example the mean commute time is n^2.4.

\subsection{Resistance of the 2D supercritical percolation cluster}

{\bf 2012 update.}
As was noted by Yoshihiro Abe, the proof in this section is wrong, and
in fact the resistance formula, applied naively, only gives that the
resistance is bounded by $C\log^2 n$. Unfortunately, we know of at
least one paper who relied on corollary 3.1 since this paper was
published. {\bf end 2012 update.}

Consider supercritical $(p > 1/2)$  bond percolation on the $n
\times n$ box of the 2D square lattice. Grimmett (private
communication) asked: show that almost surely, with respect to the
percolation measure $\P_p$, the maximal resistance between any
pair of vertices on the giant component is bounded by $C \log n$.
Denote by $R_C^n$ the maximal resistance between any pair of
vertices on the largest cluster of the percolation inside the $n
\times n$ box.  Indeed we have
\begin{coro}
$$
\P_p (R^n_C < C_p \log n) \rightarrow 1.
$$
\end{coro}

\begin{proof}
By theorem~\ref{main} it is enough to show that for $C$
sufficiently large and $c > 0 $ sufficiently small the probability
that any connected set $S$ in the giant component of size bigger
than $C \log n$, has boundary of size bigger than $ c |S|^{1/2}$
goes to $1$ with $n$. This indeed follows from an old argument of
Kesten \cite{k1} and is done explicitly in section 2.3 of
\cite{BM}.
\end{proof}

For more on the relationships between random walks and percolation clusters see
\cite{MR}, \cite{BM} and the references therein.

\section{A Conjecture}

The Cheeger constant of a finite \emph{transitive} graph is at least the
reciprocal of the diameter (see \cite{BaS}). We hope the following
stronger conjecture holds.

\begin{conj}
\label{strings} Let $G$ be finite, connected and vertex transitive.  Show that
if $\diam (G) < |G|^{\alpha}$ then $|\partial S |> c_{\alpha}
|S|^{1 - \alpha}$ for any $S$, $ 1 \leq |S| \leq |G|/2$ .

\end{conj}

If true, the first part of the next conjecture will follow along
the lines of  proof of theorem~\ref{main}.
\begin{conj}
\label{resconj}
Let $G$ be finite, connected and vertex transitive. For any
two vertices
$$
R_{v,u} < C + {\diam (G)^2 \log|G| \over |G|}.
$$
In addition, if the diameter is $o(|G|)$ then the electric
resistance between any two vertices is $o(\diam(G))$.
\end{conj}

These conjectures should be compared with the case of infinite
vertex transitive graphs which was settled by Varopoulos
\cite{VSC}, the only recurrent vertex transitive graphs are roughly
isometric to $\Z$ or $\Z^2$.

This might be the point to note that for vertex transitive graphs it is
possible to prove isoperimetric inequalities of this kind by examining balls
only, using a result of Coulhon and Saloff-Coste \cite{CSC}. Their theorem
(theorem 1 ibid.) is stated for Cayley graphs of finitely generated groups, but
it is not difficult to generalize it to edge transitive graphs (or to vertex
transitive graphs with bounded degree), e.g.\ using
theorem 6 ibid.\ with an isometry-invariant flow in the spirit of \cite{BaS}.

% I am lying here but only very slightly. If the group of isometries of the
% graph has a non-trivial point stabilizer, then it is necessary (in the
% calculations on page 13 of CSC) to multiply and divide by the size of this
% stabilizer. However, the result is still true.


\begin{thebibliography}{99}

\bibitem{AF}
D. Aldous and J. Fill,  Reversible Markov Chains and Random Walks
on Graphs.\\
\url{ http://stat-www.berkeley.edu/users/aldous/book.html}

\bibitem{BaS} L. Babai and M. Szegedy,
Local expansion of symmetrical graphs. {\it Combin. Probab.
Comput.} 1 (1992), no. 1, 1--11.

\bibitem{BM}
I. Benjamini and E. Mossel, On the mixing time of simple random
walk on the super critical percolation cluster. {\it Prob. Theor.
and Rel. Fields}, to appear.

\bibitem{CSC}
T. Coulhon and L. Saloff-Coste, Isoperimetrie pour les groupes et les varietes,
\emph{Rev. Mat. Iberoamericana} 9 (1993) 293-314. In French.

\bibitem{DS}
P. Doyle and J. Snell , Random Walks and Electric Networks.\\
\url{http://front.math.ucdavis.edu/math.PR/0001057}


\bibitem{k1}
H. Kesten, On the time constant and path length of first-passage
percolation. {\it Adv. in Appl. Probab.} 12 (1980), no. 4,
848--863.

\bibitem{LP}
R. Lyons with Y. Peres, Probability on Trees and Networks.\\
\url{http://www.math.gatech.edu/~rdlyons/prbtree/prbtree.html}

\bibitem{MR}
P. Mathieu and E. Remy, Isoperimetry and heat kernel decay on percolation
clusters, \emph{C. R. Acad. Sci. Paris S\'er. I Math.} 332 (2001), no. 10, 927--931.

\bibitem{P}
Y. Peres, Probability on Trees: An Introductory Climb.\\
\url{http://stat-www.berkeley.edu/~peres/}

\bibitem{VSC}
 N. Varopoulos,  L.  Saloff-Coste and T. Coulhon,  Analysis and geometry on groups.
 Cambridge Tracts in Mathematics, 100. Cambridge University Press, Cambridge, 1992. xii+156 pp.

\end{thebibliography}
\end{document}